# Estimation of errors on perturbation of function contractivity factors and box-counting dimension of hidden variable recurrent fractal interpolation function

Mi-Kyong Ri, Chol-Hui Yun*

Faculty of Mathematics, **Kim Il Sung** University,

Pyongyang, Democratic People's Republic of Korea

*Email Address: ch.yun@ryongnamsan.edu.kp

**Abstract:** In this paper, we study errors on perturbation of function contractivity factors and box-counting dimension of hidden variable recurrent fractal interpolation function (HVRFIF). The HVRFIF is a hidden variable fractal interpolation function (HVFIF) constructed by recurrent iterated function system (RIFS) with function contractivity factors. The contractivity factors of RIFS determine fractal characteristics and shape of its attractor, so that the HVRFIF with function contractivity factors has more flexibility and diversity than the HVFIF with constant contractivity factors. Stability of interpolation function according to perturbation of the contractivity factors and the box-counting dimension of interpolation function plays very important roles in determining whether these functions can be applied to practical problems or not. We first estimate errors on perturbation of function contractivity factors and then obtain the upper and lower bounds of the box-counting dimension of one variable HVRFIF. Finally, in the similar way, we get the lower and upper bounds of box-counting dimension of hidden variable bivariable recurrent fractal interpolation function (HVBRFIF).

**Keywords:** Recurrent iterated function system (RIFS), Hidden variable fractal interpolation function (HVFIF), Hidden variable recurrent fractal interpolation function (HVRFIF), Function contractivity factor, Box-counting dimension.

**AMS Subject Classification:** 37C45; 28A80; 41A05

## 1. Introduction

Fractal interpolation function (FIF) is an interpolation function of which graph is a fractal set. In 1986, M. F. Barnsley [1] introduced the concept of FIF to model better natural phenomena that are irregular and complicated, and many researchers have studied on construction of FIFs [3, 4, 9, 15-17, 21, 25-28 and 30], fractal dimension [6, 8, 16, 17, 19, 26 and 27] and its analytic properties [6-8, 10, 12-14, 22-25 and 29].

A recurrent fractal interpolation function (RFIF), which is defined by Barnsley ([2]), is an interpolation function of which graph is an attractor of a recurrent iterated functions system (RIFS). The RIFS is a generalization of IFS and constructs local self-similar sets that have more complicated structures than self-similar sets constructed by IFSs, and applied to image compression to improve the quality of decoded image. [5, 11 and 19] P. Bouboulis et al. [4] introduced a construction of recurrent bivariate fractal interpolation functions and estimated their box-counting dimension. Yun et al. ([26]) constructed the RFIFs using the RIFS with function vertical scaling factors and estimated the box-counting dimension of the graph of the constructed RFIFs. They also constructed bivariable RFIFs with function vertical scaling factors and estimated the box-counting dimension of the graph of the RFIFs ([27]). D.C. Luor [15] introduced a construction of FIFs with locally self-similar graphs in $R^2$, where contractivity factors are homeomorphisms that map domains to regions and are not contractions.

M.F. Barnsley *et al.* [3] introduced a concept of hidden variable fractal interpolation function (HVFIF) whose graph is neither self-similar set nor self-affine one. By changing hidden variables we can control fractal characteristics and shapes of HVFIFs more flexibly. In many articles, authors have studied HVFIFs; construction, analytic properties of HVFIF with constant contractivity factors [6, 7, 13 and 14], their box-counting dimension [6 and 19], and construction and analytic properties of HVFIF with function contractivity factors [21 and 28-30]. In order to ensure flexibility and diversity of FIFs, C.H. Yun [30] introduced a construction of one variable HVRFIF and bivariable HVRFIF using RFIS with four function contractivity factors. Furthermore analytic properties of the HVRFIFs were studied in [31].

The main aim of this paper is to estimate errors on perturbation of function contractivity factors and box-counting dimension of HVRFIF introduced in [30]. To do this, first of all, we introduce the construction of HVRFIFs constructed in [30] and estimate their errors on perturbation of function contractivity factors. Next, we obtain the upper and lower bounds of box-counting dimension of the HVRFIF. In the similar way, we get the lower and upper bounds of box-counting dimension of HVBRFIFs

The remainder of this paper is organized as follows: In section 2, we introduce the construction of HVRFIFs in [30]. (Theorem 1 and Theorem 2) In section 3, we estimate errors of HVRFIFs on perturbation of function contractivity factors. (Theorem 3, which ensures the stability of the constructed HVRFIFs according to small change of function contractivity factors) In section 4, we get the upper and lower bounds of box-counting dimension of the HVRFIF. (Theorem 5, which gives fractal structure of the HVRFIF) In section 5, we introduce the result for the box-counting dimension of HVBRFIFs (Theorem 6, which gives fractal structure of the HVBRFIF).

## 2. Hidden variable recurrent fractal interpolation function(HVRFIF)

In this section, we introduce the construction of one variable HVRFIFs in [30], that is used to estimate errors on perturbation of function contractivity factors and box-counting dimension of HVRFIFs.

### 2.1 Construction of RIFS

Let a data set $P_0$ in $R^2$ be as follows:

$$P_0 = \{(x_i, y_i) \in R^2; i = 0, 1, \cdots, n\}, \quad (-\infty < x_0 < x_1 < \cdots < x_n < +\infty).$$

We extend the dataset as follows:

$$P = \{(x_i, y_i, z_i) = (x_i, \vec{y}_i) \in R^3; i = 0, 1, \cdots, n\}, \ (-\infty < x_0 < x_1 < \cdots < x_n < \infty),$$

where $\vec{y}_i = (y_i, z_i)$ and $z_i$, $i = 0, 1, \cdots, n$ are parameters. Let denote $I_i = [x_{i-1}, x_i]$ and $I = [x_0, x_n]$ where $I_i$, $i \in \{1, \cdots, n\}$ is called a region. Let $l$ be an integer such that $2 \le l \le n$. We make subintervals $\tilde{I}_k$, $k = 1, \cdots, l$ of $I$ consisting some regions, which are called domain. Let we denote start point and end point of $\tilde{I}_k$ by $s(k)$, $e(k)$, respectively. Then we have the following mappings:

$$s : \{1, \cdots, l\} \to \{1, \cdots, n\}, \ e : \{1, \cdots, l\} \to \{1, \cdots, n\},$$

and $\tilde{I}_k$ is denoted by $\tilde{I}_k = [x_{s(k)}, x_{e(k)}]$. Suppose that $e(k) - s(k) \ge 2$, $k = 1, \cdots, l$, which means that the domain $\tilde{I}_k$ contains at least 2 of region $I_i$s. Let $f$ be a function defined on $K \subset R^n$ satisfying the Lipschitz condition. Then there is positive number $L_f$ such that for any $x, y \in K$, $\| f(x) - f(y) \| \le L_f \| x - y \|$. From now on, for all function *f* satisfying the Lipschitz condition we use notations like $L_f$. Especially, in the case of $L_f < 1$, i.e. $f$ is a contraction function, we denote $L_f$ by $c_f$.

For each $i\,(1\le i\le n)$, we take a $k(\in\{1,\cdots,l\})$ and denote it by $\gamma(i)$. Let mappings $L_{i,k}:[x_{s(k)},x_{e(k)}]\to[x_{i-1},x_i]$ be contraction homeomorphisms that map end points of $\tilde{I}_k$ to end points of $I_i$, i.e. $L_{i,k}(\{x_{s(k)},x_{e(k)}\})=\{x_{i-1},x_i\}$.

We define mappings $F_{i,k}:\tilde{I}_k\times\mathrm{R}^2\to\mathrm{R}^2,\ i=1,\cdots,n$ as follows:

$$\vec{F}_{i,k}(x,\ \vec{y})=\begin{pmatrix} s_i(L_{i,k}(x))y+s_i'(L_{i,k}(x))z+q_{i,k}(x)\\ \tilde{s}_i(L_{i,k}(x))y+\tilde{s}_i'(L_{i,k}(x))z+\tilde{q}_{i,k}(x)\end{pmatrix}=\begin{pmatrix} s_i(L_{i,k}(x)) & s_i'(L_{i,k}(x))\\ \tilde{s}_i(L_{i,k}(x)) & \tilde{s}_i'(L_{i,k}(x))\end{pmatrix}\begin{pmatrix} y\\ z\end{pmatrix}+\begin{pmatrix} q_{i,k}(x)\\ \tilde{q}_{i,k}(x)\end{pmatrix},$$

where $s_i,s_i',\tilde{s}_i,\tilde{s}_i':I_i\to R$ are Lipschitz functions on $I_i$ whose absolute value is less than 1 and $q_{i,k}$, $\tilde{q}_{i,k}:\tilde{I}_k\to R$ are Lipschitz functions such that if $\alpha\in\{s(k),e(k)\}$, $L_{i,k}(x_\alpha)=x_a$, $a\in\{i-1,i\}$, then $\vec{F}_{i,k}(x_\alpha,\ \vec{y}_\alpha)=\vec{y}_a$. Then, $\vec{F}_{i,k}(x,\ \vec{z})$ is obviously Lipschitz function.

Example 1. ([30]) An example of $q_{i,k}$ and $\tilde{q}_{i,k}$ satisfying the above conditions is as follows:

$$q_{i,k}(x)=-s_i(L_{i,k}(x))g_{i,k}(x)-s_i'(L_{i,k}(x))g_{i,k}'(x)+h_i(L_{i,k}(x)),$$

$$\tilde{q}_{i,k}(x)=-\tilde{s}_i(L_{i,k}(x))g_{i,k}(x)-\tilde{s}_i'(L_{i,k}(x))g_{i,k}'(x)+\tilde{h}_i(L_{i,k}(x)),$$

$$g_{i,k}(x)=\frac{x-x_{s(k)}}{x_{e(k)}-x_{s(k)}}y_{e(k)}+\frac{x-x_{e(k)}}{x_{s(k)}-x_{e(k)}}y_{s(k)},\quad g_{i,k}'(x)=\frac{x-x_{s(k)}}{x_{e(k)}-x_{s(k)}}z_{e(k)}+\frac{x-x_{e(k)}}{x_{s(k)}-x_{e(k)}}z_{s(k)},$$

$$h_i(x)=\frac{x-x_{i-1}}{x_i-x_{i-1}}y_i+\frac{x-x_i}{x_{i-1}-x_i}y_{i-1},\quad \tilde{h}_i(x)=\frac{x-x_{i-1}}{x_i-x_{i-1}}z_i+\frac{x-x_i}{x_{i-1}-x_i}z_{i-1}.$$

Let $D\subset R^2$ be a sufficiently large bounded set containing $\vec{y}_i$, $i=1,\cdots,n$.

We define transformations $\vec{W}_i:\tilde{I}_{\gamma(i)}\times\mathrm{D}\to I_i\times\mathrm{R}^2,\ i=1,\cdots,n$ by

$$\vec{W}_i(x,\vec{y})=(L_{i,\gamma(i)}(x),\ \vec{F}_{i,\gamma(i)}(x,\vec{y})),\quad i=1,\cdots,n.$$

Then, as we know from the definitions of $L_{i,\gamma(i)}$ and $\vec{F}_{i,\gamma(i)}$, $\vec{W}_i$ maps the data points on end of the domain $\tilde{I}_{\gamma(i)}$ to the data points of $I_i$. For a function $f$, let us denote $\bar{f}=\max_x|f(x)|$. We denote $\bar{S}=\max\{\bar{s}_i+\bar{\tilde{s}}_i,\bar{s}_i'+\bar{\tilde{s}}_i';i=1,\cdots,n\}$.

**Theorem 1([30]).** *If* $\bar{S}<1$, *then there exists some distance* $\rho_\theta$ *equivalent to the Euclidean metric such that* $\vec{W}_i\ (i=1,\cdots,n)$ *are contraction transformations with respect to the distance* $\rho_\theta$.

We define a row-stochastic matrix $M=(p_{st})_{n\times n}$ by

$$p_{st}=\begin{cases}1/a_s, & I_s\subseteq\tilde{I}_{\gamma(t)}\\ 0, & I_s\not\subset\tilde{I}_{\gamma(t)}\end{cases}\tag{1}$$

where for every $s$ ($1\le s\le n$), the number $a_s$ indicates the number of the domains $\tilde{I}_k$ containing the region $I_s$. A connection matrix $C$ is defined by

$$c_{st}=\begin{cases}1, & p_{ts}>0\\ 0, & p_{ts}=0\end{cases}.$$

It is obvious that if the row-stochastic matrix is irreducible, then the connection matrix is also irreducible.

Then, we have RIFS $\{\mathrm{R}^3;\ M;\ W_i,i=1,\cdots,n\}$ corresponding to the extended dataset *P*.

### 2.2 Construction of HVRFIF

For the RIFS constructed above, we have the following theorem.

**Theorem 2([30]).** *There is a continuous function* $\vec{f}$ *interpolating the extended data set* $P$ *such*

*that the graph of $\vec{f}$ is the attractor A of RIFS constructed above.*

The function space and the Read-Bajraktarevic operator considered in Theorem 2 are as follows:

$$\overline{C}(I)=\{\vec{h}:I\to \mathrm{R}^2;\vec{h}\ \text{ interpolates the extended data set } P \text{ and is continuous}\}$$

$$(T\vec{h})(x)=\vec{F}_{i,\gamma(i)}(L^{-1}_{i,\gamma(i)}(x),\ \vec{h}(L^{-1}_{i,\gamma(i)}(x))),\ \ x\in I_i\,.$$

Then, a unique fixed point of the Read-Bajraktarevic operator is $\vec{f}(x)=\vec{F}_{i,\gamma(i)}(L^{-1}_{i,\gamma(i)}(x),\ \vec{f}(L^{-1}_{i,\gamma(i)}(x)))$.

Let us denote the vector valued function $\vec{f}:I\to \mathrm{R}^2$ by $\vec{f}=(f_1,f_2)$, where $f_1:I\to R$ interpolates the given dataset $P_0$ and is called a hidden variable recurrent fractal interpolation function (HVRFIF), and $f_2:I\to \mathrm{R}$ interpolates the set $\{(x_i,z_i,)\in \mathrm{R}^2; i=0,1,\cdots,n\}$. Furthermore, a set $\{(x,\ f_1(x)):x\in I\}$ is a projection of $A$ on $R^2$. Then we have

$$\vec{f}(x,\ y)=\vec{F}_{i,\gamma(i)}(L^{-1}_{i,\gamma(i)}(x),\ \vec{f}(L^{-1}_{i,\gamma(i)}(x))),\ \ x\in I_i\,,$$

i.e.

$$\vec{f}(x)=\vec{F}_{i,\gamma(i)}(L^{-1}_{i,\gamma(i)}(x),\ f_1(L^{-1}_{i,\gamma(i)}(x)), f_2(L^{-1}_{i,\gamma(i)}(x))),\ \ x\in I_i\,.$$

Therefore, for every $x\in I$, the HVRFIF $f_1$ and $f_2$ satisfy the following;

$$\begin{aligned} f_1(x)&=s_i(x)f_1(L^{-1}_{i,\gamma(i)}(x))+s_i'(x)f_2(L^{-1}_{i,\gamma(i)}(x))+q_{i,\gamma(i)}(L^{-1}_{i,\gamma(i)}(x)),\\ f_2(x)&=\tilde{s}_i(x)f_1(L^{-1}_{i,\gamma(i)}(x))+\tilde{s}_i'(x)f_2(L^{-1}_{i,\gamma(i)}(x))+\tilde{q}_{i,\gamma(i)}(L^{-1}_{i,\gamma(i)}(x)),\ x\in I_i. \end{aligned} \tag{2}$$

## 3. Errors of HVRFIFs on perturbation of function contractivity factors

In this section, we estimate errors of HVRFIFs on perturbation of function contractivity factors in (2) Let $I=[0,1]$ and $L_{i,k}:\tilde{I}_k\to I_i$ be similitude mappings. We use the following notations:

$$\Omega=\max_k\{\bar{s}_k,\bar{s}_k'\},\ \ \tilde{\Omega}=\max_k\{\bar{\tilde{s}}_k,\bar{\tilde{s}}_k'\}.$$

By $q_{i,k}$, $\tilde{q}_{i,k}$ in Example 1, for any $x\in I$, we have

$$\begin{aligned} f_1(x)&=s_i(x)[f_1(L^{-1}_{i,\gamma(i)}(x))-g_{i,\gamma(i)}(L^{-1}_{i,\gamma(i)}(x))]+s_i'(x)[f_2(L^{-1}_{i,\gamma(i)}(x))-g'_{i,\gamma(i)}(L^{-1}_{i,\gamma(i)}(x))]+h_i(x)\\ f_2(x)&=\tilde{s}_i(x)[f_1(L^{-1}_{i,\gamma(i)}(x))-g_{i,\gamma(i)}(L^{-1}_{i,\gamma(i)}(x))]+\tilde{s}_i'(x)[f_2(L^{-1}_{i,\gamma(i)}(x))-g'_{i,\gamma(i)}(L^{-1}_{i,\gamma(i)}(x))]+\tilde{h}_i(x) \end{aligned} \tag{3}$$

Now we give perturbations to the function contractivity factors as follows:

$$\begin{aligned} f_1^*(x)&=(s_i(x)+\delta_i(x))[f_1^*(L^{-1}_{i,\gamma(i)}(x))-g_{i,\gamma(i)}(L^{-1}_{i,\gamma(i)}(x))]+\\ &\quad+(s_i'(x)+\delta_i'(x))[f_2^*(L^{-1}_{i,\gamma(i)}(x))-g'_{i,\gamma(i)}(L^{-1}_{i,\gamma(i)}(x))]+h_i(x),\\ f_2^*(x)&=(\tilde{s}_i(x)+\tilde{\delta}_i(x))[f_1^*(L^{-1}_{i,\gamma(i)}(x))-g_{i,\gamma(i)}(L^{-1}_{i,\gamma(i)}(x))]+\\ &\quad+(\tilde{s}_i'(x)+\tilde{\delta}_i'(x))[f_2^*(L^{-1}_{i,\gamma(i)}(x))-g'_{i,\gamma(i)}(L^{-1}_{i,\gamma(i)}(x))]+\tilde{h}_i(x) \end{aligned} \tag{4}$$

We denote $\Delta=\max_k\{\overline{\delta_k},\ \overline{\delta_k'}\}$, $\tilde{\Delta}=\max_k\{\overline{\tilde{\delta}_k},\ \overline{\tilde{\delta}_k'}\}$.

**Theorem 3.** *If* $\Omega+\Delta+\tilde{\Omega}+\tilde{\Delta}<1$, *then there exist positive numbers* $P,\ Q,\ \tilde{P}$ *and* $\tilde{Q}$ *such that*

$$|\,f_1(x)-f_1^*(x)\,|\le\frac{P\Delta+Q\tilde{\Delta}}{1-(\Omega+\Delta)-(\tilde{\Omega}+\tilde{\Delta})}, \tag{5}$$

$$|f_2(x)-f_2^*(x)| \le \frac{\tilde{P}\Delta+\tilde{Q}\tilde{\Delta}}{1-(\Omega+\Delta)-(\tilde{\Omega}+\tilde{\Delta})}, \quad x\in I \tag{6}$$

**Proof.** For any $x\in I$, there exists $x'\in \tilde{I}_{\gamma(i)}$ such that $x=L_{i,\gamma(i)}(x')$. Then we have

$$\begin{aligned}
|f_1(x)-f_1^*(x)| &= |\{s_i(x)[f_1(x')-g_{i,\gamma(i)}(x')]+s_i'(x)[f_2(x')-g_{i,\gamma(i)}'(x')]+h_i(x)\}- \\
&\quad -\{(s_i(x)+\delta_i(x))[f_1^*(x')-g_{i,\gamma(i)}(x')]+(s_i'(x)+\delta_i'(x))[f_2^*(x')-g_{i,\gamma(i)}'(x')]+h_i(x)\}| \\
&= |(s_i(x)+\delta_i(x))[f_1(x')-f_1^*(x')]+(s_i'(x)+\delta_i'(x))[f_2(x')-f_2^*(x')] \\
&\quad -\delta_i(x)[f_1(x')-g_{i,\gamma(i)}(x')]-\delta_i'(x)[f_2(x')-g_{i,\gamma(i)}'(x')]| \\
&\le |s_i(x)+\delta_i(x)|\cdot|f_1(x')-f_1^*(x')|+|s_i'(x)+\delta_i'(x)|\cdot|f_2(x')-f_2^*(x')| \\
&\quad +|\delta_i(x)|\cdot|f_1(x')-g_{i,\gamma(i)}(x')|+|\delta_i'(x)|\cdot|f_2(x')-g_{i,\gamma(i)}'(x')| \\
&\le (\Omega+\Delta)(\|f_1-f_1^*\|_\infty+\|f_2-f_2^*\|_\infty)+\Delta(\|f_1\|_\infty+\|f_2\|_\infty+\bar{g}_{i,\gamma(i)}+\bar{g}_{i,\gamma(i)}')
\end{aligned} \tag{7}$$

Now, let us estimate norm $\|\cdot\|_\infty$ in (7). From the definition of $g_{i,k}$, $g_{i,k}'$, $h_i$ and $\tilde{h}_i$ in Example 1, we have

$$\begin{aligned}
&\bar{g}_{i,k} = \max_{x\in I_k}\left\{\left|\frac{x-x_{s(k)}}{x_{e(k)}-x_{s(k)}}y_{e(k)}+\frac{x-x_{e(k)}}{x_{s(k)}-x_{e(k)}}y_{s(k)}\right|\right\} \le \max\{|y_{s(k)}|,\ |y_{e(k)}|\}, \\
&\bar{g}_{i,k}' \le \max\{|z_{s(k)}|,\ |z_{e(k)}|\} \\
&\bar{h}_i \le \max\{|y_i|,\ |y_{i-1}|\},\ \bar{\tilde{h}}_i \le \max\{|z_i|,\ |z_{i-1}|\}
\end{aligned} \tag{8}$$

and by (4), we get

$$\begin{aligned}
|f_1(x)| &\le \Omega(|f_1(x')|+|f_2(x')|)+\Omega(\bar{g}_{i,\gamma(i)}+\bar{g}_{i,\gamma(i)}')+\bar{h}_i(x) \\
&\le \Omega(\|f_1\|_\infty+\|f_2\|_\infty)+(\Omega+1)\max_i\{|y_i|\}+\Omega\max_i\{|z_i|\},
\end{aligned} \tag{9}$$

$$|f_2(x)| \le \tilde{\Omega}(\|f_1\|_\infty+\|f_2\|_\infty)+\tilde{\Omega}\max_i\{|y_i|\}+(\tilde{\Omega}+1)\max_i\{|z_i|\}. \tag{10}$$

Then by (10), we obtain

$$\|f_2\|_\infty \le \tilde{\Omega}(\|f_1\|_\infty+\|f_2\|_\infty)+\tilde{\Omega}\max_i\{|y_i|\}+(\tilde{\Omega}+1)\max_i\{|z_i|\},$$

$$\|f_2\|_\infty \le \frac{1}{1-\tilde{\Omega}}[\tilde{\Omega}\|f_1\|_\infty+\tilde{\Omega}\max_i\{|y_i|\}+(\tilde{\Omega}+1)\max_i\{|z_i|\}]$$

and by substituting it for (9), we have

$$\begin{aligned}
\|f_1\|_\infty &\le \frac{1+\Omega-\tilde{\Omega}}{1-\Omega-\tilde{\Omega}}\max_i\{|y_i|\}+\frac{2\Omega}{1-\Omega-\tilde{\Omega}}\max_i\{|z_i|\}, \\
\|f_2\|_\infty &\le \frac{2\tilde{\Omega}}{1-\Omega-\tilde{\Omega}}\max_i\{|y_i|\}+\frac{1-\Omega+\tilde{\Omega}}{1-\Omega-\tilde{\Omega}}\max_i\{|z_i|\}.
\end{aligned} \tag{11}$$

Therefore by substituting (8) and (11) for (7), we get

$$|f_1(x)-f_1^*(x)| \le (\Omega+\Delta)(\|f_1-f_1^*\|_\infty+\|f_2-f_2^*\|_\infty)+\frac{2\Delta}{1-\Omega-\tilde{\Omega}}(\max_i\{|y_i|\}+\max_i\{|z_i|\}), \tag{12}$$

$$|f_2(x)-f_2^*(x)| \le (\tilde{\Omega}+\tilde{\Delta})(\|f_1-f_1^*\|_\infty+\|f_2-f_2^*\|_\infty)+\frac{2\tilde{\Delta}}{1-\Omega-\tilde{\Omega}}(\max_i\{|y_i|\}+\max_i\{|z_i|\}). \tag{13}$$

By (13), we have

$$\|f_2-f_2^*\|_\infty \le (\tilde{\Omega}+\tilde{\Delta})(\|f_1-f_1^*\|_\infty+\|f_2-f_2^*\|_\infty)+\frac{2\tilde{\Delta}}{1-\Omega-\tilde{\Omega}}(\max_i\{|y_i|\}+\max_i\{|z_i|\}),$$

$$\| f_2 - f_2^* \|_\infty \le \frac{\tilde{\Omega}+\tilde{\Delta}}{1-\tilde{\Omega}-\tilde{\Delta}} \| f_1 - f_1^* \|_\infty + \frac{2\tilde{\Delta}}{(1-\Omega-\tilde{\Omega})(1-\tilde{\Omega}-\tilde{\Delta})} (\max_i\{|y_i|\} + \max_i\{|z_i|\}) . \tag{14}$$

Hence by substituting (14) for (12), we get the following estimation.

$$\| f_1 - f_1^* \|_\infty \le \frac{2(\Delta + \tilde{\Delta}\Omega - \Delta\tilde{\Omega})}{(1-\Omega-\tilde{\Omega})(1-\Omega-\Delta-\tilde{\Omega}-\tilde{\Delta})} (\max_i\{|y_i|\} + \max_i\{|z_i|\})$$

$$= \frac{P\Delta + Q\tilde{\Delta}}{1-(\Omega+\Delta)-(\tilde{\Omega}+\tilde{\Delta})},$$

$$\| f_2 - f_2^* \|_\infty \le \frac{\tilde{P}\Delta + \tilde{Q}\tilde{\Delta}}{1-(\Omega+\Delta)-(\tilde{\Omega}+\tilde{\Delta})},$$

where

$$P = \frac{2(1-\tilde{\Omega})}{1-\Omega-\tilde{\Omega}} (\max_i\{|y_i|\} + \max_i\{|z_i|\}), \quad Q = \frac{2\Omega}{1-\Omega-\tilde{\Omega}} (\max_i\{|y_i|\} + \max_i\{|z_i|\}),$$

$$\tilde{P} = \frac{2\tilde{\Omega}}{1-\Omega-\tilde{\Omega}} (\max_i\{|y_i|\} + \max_i\{|z_i|\}), \quad \tilde{Q} = \frac{2(1-\Omega)}{1-\Omega-\tilde{\Omega}} (\max_i\{|y_i|\} + \max_i\{|z_i|\}).$$

Therefore, we can prove (5) and (6) for any $x \in I$. □

**Remark 1.** If $s'_{ij} = \tilde{s}_{ij} = \tilde{s}'_{ij} = 0,\ \delta'_{ij} = \tilde{\delta}_{ij} = \tilde{\delta}'_{ij} = 0$ i.e.

$$f(x) = s_i(x)[f(L^{-1}_{i,\gamma(i)}(x)) - g_{i,\gamma(i)}(L^{-1}_{i,\gamma(i)}(x))] + h_i(x)$$

and

$$f^*(x) = s_i(x)[f^*(L^{-1}_{i,\gamma(i)}(x)) - g_{i,\gamma(i)}(L^{-1}_{i,\gamma(i)}(x))] + h_i(x),$$

then $| f(x) - f^*(x) | \le \dfrac{2\Delta \cdot \max_i\{|y_i|\}}{(1-\Omega)(1-\Omega-\Delta)}, \quad x \in I$.

## 4. Box-counting dimension of the HVRFIF

In this section, we get a lower and upper bounds for the box-counting dimension of the graph of the HVRFIF in the case where the data set is $P_0 = \left\{ \left( x_0 + \frac{x_n - x_0}{n} i,\, y_i \right) \in \mathrm{R}^2; i = 0, 1, \cdots, n \right\}$, the extended data set is $P = \left\{ \left( x_0 + \frac{x_n - x_0}{n} i,\, y_i,\, z_i \right) \in \mathrm{R}^3; i = 0, 1, \cdots, n \right\}$, row-stochastic matrix $M$ in (1) is irreducible and the function contractivity factors in (2) satisfies the following conditions:

$$s_i(x) s'_i(x) \ge 0, \quad \tilde{s}_i(x) \tilde{s}'_i(x) \ge 0, \quad x \in I,\ i = 1, 2, \cdots, n.$$

We denote the number of regions contained in the $k$th domain by $\eta_k$.

Let us denote the graph of the HVRFIF $f_1$ by $Gr(f_1)$. As usual, the box-counting dimension of the set $A$ by $\dim_B A$ is defined by

$$\dim_B A = \lim_{\delta \to 0} \frac{\log N_\delta(A)}{-\log \delta} \text{ (if this limit exists)},$$

where $N_\delta(A)$ is any of the following:

(i) the smallest number of closed balls of radius $\delta$ that cover the set $A$;
(ii) the smallest number of cubes of side $\delta$ that cover the set $A$;
(iii) the number of $\delta$-mesh cubes that intersect the set $A$;

(iv) the smallest number of sets of diameter at most $\delta$ that cover the set $A$;
(v) the largest number of disjoint balls of radius $\delta$ with centers in the set $A$.

**Theorem 4** (Perron-Frobenius Theorem). ([20]) *Let $A \ge 0$ be an irreducible square matrix. Then we have the following two statements.*

*(1) The spectral radius $\rho(A)$ of $A$ is an eigenvalue of $A$ and it has strictly positive eigenvector $y$ (i.e., $y_i > 0$ for all $i$).*

*(2) $\rho(A)$ increases if any element of $A$ increases.*

In association with (1) and (2), let us denote as follows:

$$\overline{\eta} = \max_k\{\eta_k\},\ \underline{\eta} = \min_k\{\eta_k\},\ \overline{\Omega}_k = \max_{x\in I_k}\{|s_k(x)|, |s_k'(x)|\},\ \overline{\tilde{\Omega}}_k = \max_{x\in I_k}\{|\tilde{s}_k(x)|, |\tilde{s}_k'(x)|\},$$

$$\underline{\Omega}_k = \min_{x\in I_k}\{|s_k(x)|, |s_k'(x)|\},\ \underline{\tilde{\Omega}}_k = \min_{x\in I_k}\{|\tilde{s}_k(x)|, |\tilde{s}_k'(x)|\},\ |I| = x_n - x_0.$$

For a function $f$ defined on $D \subset \mathrm{R}$, denote as follows:

$$R_f[D] = \sup\{|f(x_2) - f(x_1)| : x_1, x_2 \in D\}.$$

As the matrixes $(\underline{S} + \underline{\tilde{S}})C$ and $(\overline{S} + \overline{\tilde{S}})C$ are positive irreducible, by Perron-Frobenius Theorem, they have spectral radius $\underline{\rho}$ and $\overline{\rho}$ with $\underline{\rho} \le \overline{\rho}$, respectively, where $C$ is the connection matrix in (1) and

$$\underline{S} = \mathrm{diag}(\underline{\Omega}_1,\ \underline{\Omega}_2,\ \cdots,\ \underline{\Omega}_n),\ \underline{\tilde{S}} = \mathrm{diag}(\underline{\tilde{\Omega}}_1,\ \underline{\tilde{\Omega}}_2,\ \cdots,\ \underline{\tilde{\Omega}}_n),$$

$$\overline{S} = \mathrm{diag}(\overline{\Omega}_1,\ \overline{\Omega}_2,\ \cdots,\ \overline{\Omega}_n),\ \overline{\tilde{S}} = \mathrm{diag}(\overline{\tilde{\Omega}}_1,\ \overline{\tilde{\Omega}}_2,\ \cdots,\ \overline{\tilde{\Omega}}_n).$$

We have the following theorem.

**Theorem 5.** *Let $f_1(x)$ be the HVFIF in Theorem 2. Suppose that there exist three interpolation points $(x_{\alpha_1}, y_{\alpha_1})$, $(x_{\alpha_2}, y_{\alpha_2})$, $(x_{\alpha_3}, y_{\alpha_3}) \in P_0$ $(x_{\alpha_1} < x_{\alpha_2} < x_{\alpha_3})$ which are not collinear and that take $z_{\alpha_1}$, $z_{\alpha_2}$ and $z_{\alpha_3}$ such that $(y_{\alpha_i} - y_{\alpha_j})$ $(z_{\alpha_i} - z_{\alpha_j}) > 0$, $i, j = 1, 2, 3, i \ne j$ and three points $(x_{\alpha_1}, z_{\alpha_1})$, $(x_{\alpha_2}, z_{\alpha_2})$ and $(x_{\alpha_3}, z_{\alpha_3})$ are not collinear. Then the box-counting dimension of the graph of $f_1(x)$ is as follows:*

(i) *If* $\underline{\rho} > 1$, *then* $1 + \log_{\overline{\eta}} \underline{\rho} \le \dim_B Gr(f_1) \le 1 + \log_{\overline{\eta}} \overline{\rho} + (1 - \log_{\overline{\eta}} \underline{\eta})$,

(ii) *If* $\overline{\rho} \le 1$, *then* $1 \le \dim_B Gr(f_1) \le 2 - \log_{\overline{\eta}} \underline{\eta}$.

**Proof.** Firstly, we prove (i). We denote the y-axis vertical distance from the point $(x_{\alpha_2}, y_{\alpha_2})$ to the line through the points $(x_{\alpha_1}, y_{\alpha_1})$ and $(x_{\alpha_3}, y_{\alpha_3})$ by $H$ and the $z$-axis vertical distance from the point $(x_{\alpha_2}, z_{\alpha_2})$ to the line through the points $(x_{\alpha_1}, z_{\alpha_1})$ and $(x_{\alpha_3}, z_{\alpha_3})$ by $\tilde{H}$. Then obviously $H \cdot \tilde{H} > 0$.

We apply $W_i$ to interval $I$ one time. Then we have

$$R_{f_1}[I_i] \le \overline{\Omega}_i (R_{f_1}[\tilde{I}_{\gamma(i)}] + R_{f_2}[\tilde{I}_{\gamma(i)}]) + (L_{s_i} \cdot \overline{f_1} + L_{s_i'} \cdot \overline{f_2} + \eta_i L_{q_i}) \cdot \frac{|I|}{n}$$

$$= \overline{\Omega}_i (R_{f_1}[\tilde{I}_{\gamma(i)}] + R_{f_2}[\tilde{I}_{\gamma(i)}]) + M_i \frac{|I|}{n} \ , \quad (15)$$

$$R_{f_2}[I_i] \le \overline{\tilde{\Omega}}_i (R_{f_1}[\tilde{I}_{\gamma(i)}] + R_{f_2}[\tilde{I}_{\gamma(i)}]) + \tilde{M}_i \frac{|I|}{n},$$

where $M_i = L_{s_i} \cdot \overline{f_1} + L_{s_i'} \cdot \overline{f_2} + \eta_i L_{q_i}$, $\tilde{M}_i = L_{\tilde{s}_i} \cdot \overline{f_1} + L_{\tilde{s}_i'} \cdot \overline{f_2} + \eta_i L_{\tilde{q}_i}$ .

Since $Gr(f_1)$ is the graph of a continuous function defined on $I$ , the smallest number of $\varepsilon_r$-mesh squares that cover $I_i \times \mathrm{R} \cap Gr(f_1)$ is greater than the smallest number of $\varepsilon_r$- mesh squares necessary to cover the vertical line whose length is $\underline{\Omega}_i (H + \tilde{H})$ and is less than the smallest number of $\varepsilon_r$-mesh squares necessary to cover the rectangle $I_i \times R_{f_1}[I_i]$.

By (15), we have

$$\sum_{i=1}^{n} \left\lceil \frac{\underline{\Omega}_i (H+\tilde{H})}{\varepsilon_r} \right\rceil \le N(\varepsilon_r) \le \sum_{i=1}^{n} \left( \left\lceil \frac{\overline{\Omega}_i (R_{f_1}[\tilde{I}_{\gamma(i)}] + R_{f_2}[\tilde{I}_{\gamma(i)}]) + M_i \cdot |I|/n}{\varepsilon_r} \right\rceil + 1 \right) \left( \left\lceil \frac{|I|}{n\varepsilon_r} \right\rceil + 1 \right),$$

$$\sum_{i=1}^{n} \left( \frac{\underline{\Omega}_i (H+\tilde{H})}{\varepsilon_r} - 1 \right) \le N(\varepsilon_r) \le \sum_{i=1}^{n} \left( \frac{\overline{\Omega}_i (R_{f_1}[\tilde{I}_{\gamma(i)}] + R_{f_2}[\tilde{I}_{\gamma(i)}]) + M_i \cdot |I|/n}{\varepsilon_r} + 1 \right) \left( \frac{|I|}{n\varepsilon_r} + 1 \right),$$

$$\frac{\Phi(\mathrm{H}_1)}{\varepsilon_r} - n \le N(\varepsilon_r) \le \left( \frac{\Phi(\mathrm{U}_1)}{\varepsilon_r} + n \right) \left( \frac{|I|}{n\varepsilon_r} + 1 \right),$$

where $\Phi(a) = \sum_{i=1}^{n} a_i \ \ (a = (a_1, a_2, \cdots, a_n)^T)$ and

$$\mathrm{H}_1 = (\underline{\Omega}_1 (H+\tilde{H}), \underline{\Omega}_2 (H+\tilde{H}), \cdots, \underline{\Omega}_n (H+\tilde{H}))^T, \ \ \tilde{\mathrm{H}}_1 = (\underline{\tilde{\Omega}}_1 (H+\tilde{H}), \underline{\tilde{\Omega}}_2 (H+\tilde{H}), \cdots, \underline{\tilde{\Omega}}_n (H+\tilde{H}))^T,$$

$$\mathrm{U}_1 = \begin{pmatrix} \overline{\Omega}_1 (R_{f_1}[\tilde{I}_{\gamma(1)}] + R_{f_2}[\tilde{I}_{\gamma(1)}]) \\ \overline{\Omega}_2 (R_{f_1}[\tilde{I}_{\gamma(2)}] + R_{f_2}[\tilde{I}_{\gamma(2)}]) \\ \cdots \\ \overline{\Omega}_n (R_{f_1}[\tilde{I}_{\gamma(n)}] + R_{f_2}[\tilde{I}_{\gamma(n)}]) \end{pmatrix} + \frac{|I|}{n} \begin{pmatrix} M_1 \\ M_2 \\ \cdots \\ M_n \end{pmatrix}, \ \ \tilde{\mathrm{U}}_1 = \begin{pmatrix} \overline{\tilde{\Omega}}_1 (R_{f_1}[\tilde{I}_{\gamma(1)}] + R_{f_2}[\tilde{I}_{\gamma(1)}]) \\ \overline{\tilde{\Omega}}_2 (R_{f_1}[\tilde{I}_{\gamma(2)}] + R_{f_2}[\tilde{I}_{\gamma(2)}]) \\ \cdots \\ \overline{\tilde{\Omega}}_n (R_{f_1}[\tilde{I}_{\gamma(n)}] + R_{f_2}[\tilde{I}_{\gamma(n)}]) \end{pmatrix} + \frac{|I|}{n} \begin{pmatrix} \tilde{M}_1 \\ \tilde{M}_2 \\ \cdots \\ \tilde{M}_n \end{pmatrix}.$$

Let denote that $M = (M_1, \ M_2, \cdots, \ M_n)^T$, $\tilde{M} = (\tilde{M}_1, \ \tilde{M}_2, \cdots, \ \tilde{M}_n)^T$.

Let us apply $W_j$ to every subinterval $I_i$ once again. Then we have

$$R_{f_1}[I_{ij}] \le \overline{\Omega}_j (R_{f_1}[\tilde{I}_{\gamma(j)} \cap I_i] + R_{f_2}[\tilde{I}_{\gamma(j)} \cap I_i]) + \frac{|I|}{n\eta_{\gamma(j)}} M,$$

$$R_{f_2}[I_{ij}] \le \overline{\tilde{\Omega}}_j (R_{f_1}[\tilde{I}_{\gamma(j)} \cap I_i] + R_{f_2}[\tilde{I}_{\gamma(j)} \cap I_i]) + \frac{|I|}{n\eta_{\gamma(j)}} \tilde{M}. \quad (16)$$

Since the smallest number of $\varepsilon_r$-mesh squares that cover $I_{ij} \times \mathrm{R} \cap Gr(f_1)$ is greater than the smallest number of $\varepsilon_r$- mesh squares necessary to cover the vertical line whose length is $\underline{\Omega}_j (\underline{\Omega}_i + \underline{\tilde{\Omega}}_i)(H + \tilde{H})$ and less than the smallest number of $\varepsilon_r$-mesh squares necessary to cover the rectangle $I_{ij} \times R_{f_1}[I_{ij}]$, by (16), we have

$$\sum_{i=1}^{n} \sum_{j=l'}^{l'+\eta_{\gamma(j)}-1} \left\lceil \frac{\underline{\Omega}_j (\underline{\Omega}_i + \underline{\tilde{\Omega}}_i)(H+\tilde{H})}{\varepsilon_r} \right\rceil \le N(\varepsilon_r) \le \sum_{i=1}^{n} \sum_{j=l'}^{l'+\eta_{\gamma(j)}-1} \left( \left\lceil \frac{R_{f_1}[I_{ji}]}{\varepsilon_r} \right\rceil + 1 \right) \left( \left\lceil \frac{|I|}{n\eta_{\gamma(j)}\varepsilon_r} \right\rceil + 1 \right)$$

$$\sum_{i=1}^{n} \left( \sum_{j=1}^{n} \frac{C_{ij} \underline{\Omega}_j (\underline{\Omega}_i + \underline{\tilde{\Omega}}_i)(H+\tilde{H})}{\varepsilon_r} - \eta_{\gamma(j)} \right) \le N(\varepsilon_r) \le$$

$$\le \sum_{i=1}^{n} \left( \sum_{j=1}^{n} \frac{C_{ij} [\overline{\Omega}_j (R_{f_1}[I_i] + R_{f_2}[I_i]) + M|I|/(n\eta_{\gamma(j)})]}{\varepsilon_r} + \eta_{\gamma(j)} \right) \left( \frac{|I|}{n\eta_{\gamma(j)}\varepsilon_r} + 1 \right),$$

$$\frac{\Phi(\mathrm{H}_2)}{\varepsilon_r} - n\overline{\eta} \le N(\varepsilon_r) \le \left(\frac{\Phi(\mathrm{U}_2)}{\varepsilon_r} + n\overline{\eta}\right)\left(\frac{|I|}{n\underline{\eta}\varepsilon_r} + 1\right),$$

where $\mathrm{H}_2 = \underline{S}C(\mathrm{H}_1 + \tilde{\mathrm{H}}_1)$, $\tilde{\mathrm{H}}_2 = \underline{\tilde{S}}C(\mathrm{H}_1 + \tilde{\mathrm{H}}_1)$, $\mathrm{U}_2 = \overline{S}C(\mathrm{U}_1 + \tilde{\mathrm{U}}_1) + \frac{|I|}{n}M$ and $\mathrm{U}_2 = \overline{\tilde{S}}C(\mathrm{U}_1 + \tilde{\mathrm{U}}_1) + \frac{|I|}{n}\tilde{M}$.

Suppose that

$$\varepsilon_r < \frac{|I|}{n\overline{\eta}^{k-1}} \le \overline{\eta}\varepsilon_r. \tag{17}$$

If we apply $W_i$ to $I$ $k$ times, we have

$$\frac{\Phi(\mathrm{H}_k)}{\varepsilon_r} - n\overline{\eta}^{k-1} \le N(\varepsilon_r) \le \left(\frac{\Phi(\mathrm{U}_k)}{\varepsilon_r} + n\overline{\eta}^{k-1}\right)\left[\frac{|I|}{n\overline{\eta}^{k-1}\varepsilon_r}\cdot\left(\frac{\overline{\eta}}{\underline{\eta}}\right)^{k-1} + 1\right], \tag{18}$$

where

$$\mathrm{H}_k = \underline{S}C(\mathrm{H}_{k-1} + \tilde{\mathrm{H}}_{k-1}),\ \tilde{\mathrm{H}}_k = \underline{\tilde{S}}C(\mathrm{H}_{k-1} + \tilde{\mathrm{H}}_{k-1}),$$

$$\mathrm{U}_k = \overline{S}C(\mathrm{U}_{k-1} + \tilde{\mathrm{U}}_{k-1}) + \frac{|I|}{n}M,\ \mathrm{U}_k = \overline{\tilde{S}}C(\mathrm{U}_{k-1} + \tilde{\mathrm{U}}_{k-1}) + \frac{|I|}{n}\tilde{M}.$$

Then we know that

$$\mathrm{H}_k = \underline{S}C(\mathrm{H}_{k-1} + \tilde{\mathrm{H}}_{k-1}) = \underline{S}C(\underline{S}C + \underline{\tilde{S}}C)(\mathrm{H}_{k-2} + \tilde{\mathrm{H}}_{k-2}) = \underline{S}C(\underline{S}C + \underline{\tilde{S}}C)^{k-2}(\mathrm{H}_1 + \tilde{\mathrm{H}}_1), \tag{19}$$

$$\mathrm{U}_k = \overline{S}C(\mathrm{U}_{k-1} + \tilde{\mathrm{U}}_{k-1}) + \frac{|I|}{n}M = \overline{S}C(\overline{S}C + \overline{\tilde{S}}C)(\mathrm{U}_{k-2} + \tilde{\mathrm{U}}_{k-2}) + \frac{|I|}{n}\overline{S}C(M + \tilde{M}) + \frac{|I|}{n}M = \cdots =$$

$$= \overline{S}C(\overline{S}C + \overline{\tilde{S}}C)^{k-2}(\mathrm{U}_1 + \tilde{\mathrm{U}}_1) + \frac{|I|}{n}\overline{S}C(\overline{S}C + \overline{\tilde{S}}C)^{k-3}(M + \tilde{M}) + \cdots + \frac{|I|}{n}\overline{S}C(M + \tilde{M}) + \frac{|I|}{n}M. \tag{20}$$

By Perron-Frobenius theorem, the matrixes $(\underline{S} + \underline{\tilde{S}})C$ and $(\overline{S} + \overline{\tilde{S}})C$ have strictly positive eigenvectors $\underline{\mathrm{e}}$ and $\overline{\mathrm{e}}$ which correspond to spectral radius $\underline{\rho}$ and $\overline{\rho}$, respectively, such that

$$\underline{\mathrm{e}} \le \mathrm{H}_1 + \tilde{\mathrm{H}}_1,\ \overline{\mathrm{e}} \ge \mathrm{U}_1 + \tilde{\mathrm{U}}_1,\ \overline{\mathrm{e}} \ge \frac{|I|}{n}(M + \tilde{M}).$$

Then by (18)-(20), we get

$$\frac{\underline{\rho}^{k-2}\Phi(\underline{S}C\underline{\mathrm{e}})}{\varepsilon_r} - n\overline{\eta}^{k-1} \le N(\varepsilon_r) \le \left(\frac{\Phi(\overline{\mathrm{e}})(\overline{\rho}^{k-1} + \cdots + \overline{\rho} + 1)}{\varepsilon_r} + n\overline{\eta}^{k-1}\right)\left[\frac{|I|}{n\overline{\eta}^{k-1}\varepsilon_r}\cdot\left(\frac{\overline{\eta}}{\underline{\eta}}\right)^{k-1} + 1\right] \tag{21}$$

By the assumption, since $1 < \underline{\rho} \le \overline{\rho}$, we can get

$$\frac{\underline{\rho}^{k-2}\Phi(\underline{S}C\underline{\mathrm{e}})}{\varepsilon_r} - n\overline{\eta}^{k-1} \le N(\varepsilon_r) \le \left(\frac{\Phi(\overline{\mathrm{e}})}{\varepsilon_r}\frac{\overline{\rho}^{k-1}(1 - 1/\overline{\rho}^k)}{1 - 1/\overline{\rho}} + n\overline{\eta}^{k-1}\right)\left[\frac{|I|}{n\overline{\eta}^{k-1}\varepsilon_r}\cdot\left(\frac{\overline{\eta}}{\underline{\eta}}\right)^{k-1} + 1\right],$$

and by (17), we have

$$\frac{\underline{\rho}^{k-2}}{\varepsilon_r}\left(\Phi(\underline{S}C\underline{\mathrm{e}}) - \frac{|I|}{\underline{\rho}^{k-2}}\right) \le N(\varepsilon_r) \le \frac{\overline{\rho}^{k-1}}{\varepsilon_r}\left(\frac{\overline{\eta}}{\underline{\eta}}\right)^{k-1}\left(\Phi(\overline{\mathrm{e}})\frac{1 - 1/\overline{\rho}^k}{1 - 1/\overline{\rho}} + \frac{|I|}{\overline{\rho}^{k-1}}\right)\left[\overline{\eta} + \left(\frac{\underline{\eta}}{\overline{\eta}}\right)^{k-1}\right],$$

$$\frac{\underline{\rho}^{k-2}}{\varepsilon_r}\alpha \le N(\varepsilon_r) \le \frac{\overline{\rho}^{k-1}}{\varepsilon_r}\left(\frac{\overline{\eta}}{\underline{\eta}}\right)^{k-1}\beta$$

where $\alpha = \Phi(\underline{S}C\underline{\mathrm{e}}) - \frac{|I|}{\underline{\rho}^{k-2}}$ and $\beta = \left(\Phi(\overline{\mathrm{e}})\frac{1 - 1/\overline{\rho}^k}{1 - 1/\overline{\rho}} + \frac{|I|}{\overline{\rho}^{k-1}}\right)\left[\overline{\eta} + \left(\frac{\underline{\eta}}{\overline{\eta}}\right)^{k-1}\right] > 0$. Hence we have

$$1 - \frac{(k-2)\log\underline{\rho}}{\log\varepsilon_r} - \frac{\log\alpha}{\log\varepsilon_r} \le \frac{\log N(\varepsilon_r)}{-\log\varepsilon_r} \le 1 - \frac{(k-1)\log\overline{\rho}}{\log\varepsilon_r} - \frac{(k-1)\log(\overline{\eta}/\underline{\eta})}{\log\varepsilon_r} - \frac{\log\beta}{\log\varepsilon_r},$$

$$1 + \frac{(k-2)\log\underline{\rho}}{(k-2)\log\overline{\eta} + \log n} - \frac{\log\alpha}{\log\varepsilon_r} \le \frac{\log N(\varepsilon_r)}{-\log\varepsilon_r} \le 1 + \frac{(k-1)\log\overline{\rho}}{(k-1)\log\overline{\eta} + \log n} + \frac{(k-1)\log(\overline{\eta}/\underline{\eta})}{(k-1)\log\overline{\eta} + \log n} - \frac{\log\beta}{\log\varepsilon_r},$$

$$1 + \log_{\overline{\eta}}\underline{\rho} \le \lim_{\varepsilon_r \to 0}\frac{\log N(\varepsilon_r)}{-\log\varepsilon_r} \le 1 + \log_{\overline{\eta}}\overline{\rho} + (1 - \log_{\overline{\eta}}\underline{\eta}).$$

As a result, we prove (i).

Secondly, we prove (ii). Since $\overline{\rho} \le 1$, by (21), we have

$$N(\varepsilon_r) \le \left( \frac{\Phi(\mathrm{e})}{\varepsilon_r} \frac{1-\overline{\rho}^k}{1-\overline{\rho}} + n\overline{\eta}^{k-1} \right) \left[ \frac{|I|}{n\overline{\eta}^{k-1}\varepsilon_r} \cdot \left( \frac{\overline{\eta}}{\underline{\eta}} \right)^{k-1} + 1 \right],$$

$$N(\varepsilon_r) \le \frac{1}{\varepsilon_r} \left( \frac{\overline{\eta}}{\underline{\eta}} \right)^{k-1} \left( \Phi(\mathrm{e}) \frac{1-\overline{\rho}^k}{1-\overline{\rho}} + 1 \right) \left[ \overline{\eta} + \left( \frac{\eta}{\overline{\overline{\eta}}} \right)^{k-1} \right] = \frac{1}{\varepsilon_r} \left( \frac{\overline{\eta}}{\underline{\eta}} \right)^{k-1} \gamma,$$

where $\gamma = \left( \Phi(\mathrm{e}) \frac{1-\overline{\rho}^k}{1-\overline{\rho}} + 1 \right) \left[ \overline{\eta} + \left( \frac{\eta}{\overline{\overline{\eta}}} \right)^{k-1} \right] > 0$.

Then we get

$$\frac{\log N(\varepsilon_r)}{-\log \varepsilon_r} \le 1 - \frac{(k-1)\log(\overline{\eta}/\underline{\eta})}{\log \varepsilon_r} - \frac{\log \gamma}{\log \varepsilon_r} \le 1 + \frac{(k-1)\log(\overline{\eta}/\underline{\eta})}{(k-1)\log \overline{\eta} + \log n} - \frac{\log \gamma}{\log \varepsilon_r}$$

and

$$\lim_{\varepsilon_r \to 0} \frac{\log N(\varepsilon_r)}{-\log \varepsilon_r} \le 1 + \frac{\log(\overline{\eta}/\underline{\eta})}{\log \overline{\eta}} = 2 - \log_{\overline{\eta}} \underline{\eta}. \quad \square$$

**Remark 2.** If $\overline{\eta} = \underline{\eta} = \eta$, then we get the following result:

**(i)** ***If*** $\underline{\rho} > 1$**,** ***then*** $1 + \log_\eta \underline{\rho} \le \dim_B Gr(f_1) \le 1 + \log_\eta \overline{\rho}$**,**

**(ii)** ***If*** $\overline{\rho} \le 1$**,** ***then*** $\dim_B Gr(f_1) = 1$**.**

## 5. Box-counting dimension of the HVBRFIF

The method of estimating the box-counting dimension of HVBRFIFs is similar to one in Theorem 5. Therefore, we present only the result for the box-counting dimension of HVBRFIFs.

5.1. Construction of HVBRFIFs. (see [30])

Let a dataset $P_0$ on rectangular grids be given as follows:

$$P_0 = \{(x_i,\ y_j,\ z_{ij}) \in \mathrm{R}^3;\ i = 0,\ 1, \cdots,\ n,\ j = 0,\ 1,\ \cdots,\ m\},\ (x_0 < x_1 < \cdots < x_n,\ y_0 < y_1 < \cdots < y_m)$$

We extend the dataset to the following one:

$$P = \{(x_i, y_j, z_{ij}, t_{ij}) = (\vec{x}_{ij}, \vec{z}_{ij}) \in \mathbf{R}^4; i = 0, 1, \cdots, n,\ j = 0, 1, \cdots, m\}, (x_0 < x_1 < \cdots < x_n,\ y_0 < y_1 < \cdots < y_m)$$

where $\vec{x}_{ij} = (x_i, y_j)$, $\vec{z}_{ij} = (z_{ij}, t_{ij})$ and $t_{ij}$, $i = 0, 1, \cdots, n,\ j = 0, 1, \cdots, m$ are parameters. We denote $N = n \cdot m$, $I_{x_i} = [x_{i-1},\ x_i]$, $I_{y_j} = [y_{j-1},\ y_j]$, $N_{nm} = \{1, \cdots, n\} \times \{1, \cdots, m\}$, $I_x = [x_0,\ x_n]$, $I_y = [y_0,\ y_m]$, $E = I_x \times I_y$, $E_{ij} = I_{x_i} \times I_{y_j}$, where $E_{ij}$ is a region. Let $l$ be an integer with $2 \le l \le N$. Next, we take rectangulars $\tilde{E}_k,\ k = 1, \cdots, l$ consisting of some regions from $E$. $\tilde{E}_k$ is a domain. Then we have $\tilde{E}_k = \tilde{I}_{x,k} \times \tilde{I}_{y,k}$, where $\tilde{I}_{x,k}$, $\tilde{I}_{y,k}$ are closed intervals on *x*-axis and *y*-axis, respectively. Since the endpoints of $\tilde{I}_{x,k}$ $(k \in \{1, \cdots, l\})$ are coincided with some endpoints of $I_{x_i}$, $i = 1, \cdots, n$, denoting start point and end point of $\tilde{I}_{x,k}$ by $s_x(k),\ e_x(k)$, respectively, we can define the mapping $s_x : \{1, \cdots, l\} \to \{1, \cdots, n\}$, $e_x : \{1, \cdots, l\} \to \{1, \cdots, n\}$. Similarly, for $\tilde{I}_{y,k}$, we define the mapping $s_y : \{1, \cdots, l\} \to \{1, \cdots, m\}$, $e_y : \{1, \cdots, l\} \to \{1, \cdots, m\}$. Then we have $\tilde{I}_{x,k} = [x_{s_x(k)}, x_{e_x(k)}]$, $\tilde{I}_{y,k} = [y_{s_y(k)}, y_{e_y(k)}]$, where we assume that $e_x(k) - s_x(k) \ge 2$, $e_y(k) - s_y(k) \ge 2$, $k = 1, \cdots, l$ which means that $\tilde{I}_{x,k}$, $\tilde{I}_{y,k}$ are intervals containing more than 2 $I_{x_i}s$, $I_{y_j}s$, respectively.

For $(i,\ j) \in N_{nm}$, we take $k (\in \{1, \cdots, l\})$ and denote it by $\gamma(i,\ j)$.

We define mappings $L_{x_i,k} : [x_{s_x(k)}, x_{e_x(k)}] \to [x_{i-1}, x_i]$, $L_{y_j,k} : [y_{s_y(k)}, y_{e_y(k)}] \to [y_{j-1}, y_j]$, $(i,\ j) \in N_{nm}$ as contraction homeomorphisms that map end points of $\tilde{I}_{x,k}$, $\tilde{I}_{y,k}$ to end points of $I_{x_i}$, $I_{y_j}$, i.e.

$$L_{x_i,k}(\{x_{s_x(k)}, x_{e_x(k)}\}) = \{x_{i-1}, x_i\}, \quad L_{y_j,k}(\{y_{s_y(k)}, y_{e_y(k)}\}) = \{y_{j-1}, y_j\}.$$

Next, we define transformations $\vec{L}_{ij,k} : \tilde{E}_k \to E_{ij}$ by $\vec{L}_{ij,k}(\vec{x}) = (L_{x_i,k}(x), L_{y_j,k}(y))$. Then, for $\alpha \in \{s_x(k), e_x(k)\}$, $\beta \in \{s_y(k), e_y(k)\}$, $\vec{L}_{ij,k}(\vec{x}_{\alpha\beta}) = \vec{x}_{ab}$ ($a \in \{i-1, i\}$ , $b \in \{j-1, j\}$).

We define mappings $\vec{F}_{ij,k} : \tilde{E}_k \times \mathrm{R}^2 \to \mathrm{R}^2$, $i = 1, \cdots, n$, $j = 1, \cdots, m$ by

$$\vec{F}_{ij,k}(\vec{x}, \vec{z}) = \begin{pmatrix} s_{ij}(\vec{L}_{ij,k}(\vec{x}))z + s'_{ij}(\vec{L}_{ij,k}(\vec{x}))t + q_{ij,k}(\vec{x}) \\ \tilde{s}_{ij}(\vec{L}_{ij,k}(\vec{x}))z + \tilde{s}'_{ij}(\vec{L}_{ij,k}(\vec{x}))t + \tilde{q}_{ij,k}(\vec{x}) \end{pmatrix} = \begin{pmatrix} s_{ij}(\vec{L}_{ij,k}(\vec{x})) & s'_{ij}(\vec{L}_{ij,k}(\vec{x})) \\ \tilde{s}_{ij}(\vec{L}_{ij,k}(\vec{x})) & \tilde{s}'_{ij}(\vec{L}_{ij,k}(\vec{x})) \end{pmatrix} \begin{pmatrix} z \\ t \end{pmatrix} + \begin{pmatrix} q_{ij,k}(\vec{x}) \\ \tilde{q}_{ij,k}(\vec{x}) \end{pmatrix},$$

where $s_{ij}, s'_{ij}, \tilde{s}_{ij}, \tilde{s}'_{ij} : E_{ij} \to \mathrm{R}$ are arbitrary Lipschitz functions whose absolute values are less than 1 and $q_{ij,k}$, $\tilde{q}_{ij,k} : \tilde{E}_k \to \mathrm{R}$ are defined as Lipschitz functions satisfying the following condition: for $\alpha \in \{s_x(k), e_x(k)\}$, $\beta \in \{s_y(k), e_y(k)\}$, $a \in \{i-1, i\}$, $b \in \{j-1, j\}$, $L_{x_i,k}(x_\alpha) = x_a$, $L_{y_j,k}(y_\beta) = y_b$, $\vec{F}_{ij,k}(\vec{x}_{\alpha\beta}, \vec{z}_{\alpha\beta}) = \vec{z}_{ab}$.

Then, $\vec{F}_{ij,k}(\vec{x}, \vec{z})$ are Lipschitz mappings.

We denote $\vec{F}_{ij,k}(\vec{x}, \vec{z})$ by $\vec{F}_{ij,k}(\vec{x}, \vec{z}) = \vec{S}_{ij,k}\vec{z} + \vec{Q}_{ij,k}(\vec{x})$, where

$$\vec{S}_{ij,k}(\vec{x}) = \begin{pmatrix} s_{ij}(\vec{L}_{ij,k}(\vec{x})) & s'_{ij}(\vec{L}_{ij,k}(\vec{x})) \\ \tilde{s}_{ij}(\vec{L}_{ij,k}(\vec{x})) & \tilde{s}'_{ij}(\vec{L}_{ij,k}(\vec{x})) \end{pmatrix}, \quad \vec{Q}_{ij}(\vec{x}) = \begin{pmatrix} q_{ij,k}(\vec{x}) \\ \tilde{q}_{ij,k}(\vec{x}) \end{pmatrix}, \quad \vec{z} = \begin{pmatrix} z \\ t \end{pmatrix}.$$

In the future, we denote simply $\vec{F}_{ij,k}$, $\vec{L}_{ij,k}, \vec{S}_{ij,k}, \vec{Q}_{ij,k}$ by $\vec{F}_{ij}$, $\vec{L}_{ij}$, $\vec{S}_{ij}, \vec{Q}_{ij}$, respectively. Let $\mathrm{D} \subset \mathrm{R}^2$ be a sufficiently large bounded set containing $\vec{z}_{ij}$, $i = 1, \cdots, n$, $j = 1, \cdots, m$.

Now, we define transformations $\vec{W}_{ij} : \tilde{E}_k \times \mathrm{D} \to E_{ij} \times \mathrm{R}^2$, $i = 1, \cdots, n$, $j = 1, \cdots, m$ by

$$\vec{W}_{ij}(\vec{x}, \vec{z}) = (\vec{L}_{ij}(\vec{x}), \ \vec{F}_{ij}(\vec{x}, \vec{z})), \quad i = 1, \cdots, n, \ j = 1, \cdots, m$$

We denote $\bar{S} = \max\{\bar{s}_{ij} + \bar{\tilde{s}}_{ij}, \bar{s}'_{ij} + \bar{\tilde{s}}'_{ij}; i = 1, \cdots, n, j = 1, \cdots, m\}$. If $\bar{S} < 1$, then there exists some distance $\rho_\theta$ equivalent to the Euclidean metric such that $\vec{W}_{ij}$, $i = 1, \cdots, n$, $j = 1, \cdots, m$ are contraction transformations with respect to the distance $\rho_\theta$. (See Theorem 3 in [30])

We define a row-stochastic matrix $M = (p_{st})_{N \times N}$ by

$$p_{st} = \begin{cases} 1/b_s, & E_{\tau^{-1}(s)} \subseteq \tilde{E}_{\gamma(\tau^{-1}(t))} \\ 0, & E_{\tau^{-1}(s)} \not\subset \tilde{E}_{\gamma(\tau^{-1}(t))} \end{cases}, \tag{22}$$

$\tau : N_{nm} \to \{1, \cdots, N\}$ is an one to one mapping defined by $\tau(i, j) = i + (j-1)n$ and the number $b_s$ indicates the number of the domains $\tilde{E}_k$, $k = 1, \cdots, l$ containing the region $E_{\tau^{-1}(s)}$, which means that $p_{st}$ is positive if there is a transformation $\vec{W}_{ij}$ mapping $E_s$ to $E_t$. The connection matrix of the row-stochastic matrix is the same as one in (1).

Then, RIFS $\{\mathrm{R}^4; M; \vec{W}_{ij}, i = 1, \cdots, n, j = 1, \cdots, m\}$ is a RIFS corresponding to the extended dataset $P$. An attractor of the RIFS is called a recurrent fractal set.

We present a sufficient condition for a fixed point of the Read-Bajraktarevic operator defined by the RIFS to interpolate the extended dataset $P$ and have a graph which is the attractor of the RIFS .

In the mapping $\vec{F}_{ij}(\vec{x}, \vec{z}) = \vec{S}_{ij}\vec{z} + \vec{Q}_{ij}(\vec{x})$ defined above, we define $\vec{Q}_{ij}$ as one satisfying condition that for some continuous function $\vec{g}$ interpolating the data set $P$, i.e. $\vec{g}(\vec{x}_{ij}) = \vec{z}_{ij}$, $i = 0, 1, \cdots, n$, $j = 0, 1, \cdots, m$,

$$\vec{F}_{ij}(x_\alpha, y, \vec{g}(x_\alpha, y)) = \vec{g}(\vec{L}_{ij}(x_\alpha, y)), \ \alpha \in \{s_x(k), e_x(k)\}, \tag{23}$$

$$\vec{F}_{ij}(x, y_\beta, \vec{g}(x, y_\beta)) = \vec{g}(\vec{L}_{ij}(x, y_\beta)), \ \beta \in \{s_y(k), e_y(k)\}. \tag{24}$$

One example is as follows:

$$\vec{Q}_{ij}(\vec{x}) = -\vec{S}_{ij}(\vec{L}_{ij}(\vec{x}))\vec{l}_{ij}(\vec{x}) + \vec{r}_{ij}(\vec{L}_{ij}(\vec{x})),$$

where $\vec{r}_{ij}(\vec{x})$ and $\vec{l}_{ij}(\vec{x})$ are Lipschitz functions satisfying the following conditions:

$$\vec{r}_{ij}(\vec{x}) = \vec{g}(\vec{x}),\ \vec{x} \in \partial E_{ij},$$
$$\vec{l}_{ij}(\vec{x}) = \vec{g}(\vec{x}),\ \vec{x} \in \partial \tilde{E}_k.$$

and we have $\vec{F}_{ij}(\vec{x}, \vec{z}) = \vec{S}_{ij}(\vec{L}_{ij}(\vec{x}))(\vec{z} - \vec{l}_{ij}(\vec{x})) + \vec{r}_{ij}(\vec{L}_{ij}(\vec{x}))$. Let us denote the attractor of the RIFS satis fying conditions (23) and (24) by **B**. Then, there is a continuous function $\vec{f}$ which interpolates the dat a set P and whose graph is the attractor **B**. (See Theorem 4 in [30])

The vector valued function $\vec{f}: E \to \mathrm{R}^2$ in Theorem 4 is denoted by $\vec{f} = (f_1(\vec{x}), f_2(\vec{x}))$, where $f_1: E \to \mathrm{R}$ interpolates the given dataset $P_0$ and is called hidden variable bivariable recurrent fractal interpolation function (HVBRFIF) for the data set $P_0$. $f_2(\vec{x})$ interpolates the set $\{(x_i, y_j, t_{ij}) = (\vec{x}_{ij}, t_{ij}) \in \mathrm{R}^3;$ $i = 0, 1, \cdots, n,\ j = 0, 1, \cdots, m\}$.

We get $\vec{f}(x, y) = \vec{F}_{ij}(\vec{L}_{ij}^{-1}(x, y),\ \vec{f}(\vec{L}_{ij}^{-1}(x, y))),\ (x, y) \in E_{ij}$, i.e.

$$\vec{f}(x, y) = \vec{F}_{ij}(\vec{L}_{ij}^{-1}(x, y),\ f_1(\vec{L}_{ij}^{-1}(x, y)), f_2(\vec{L}_{ij}^{-1}(x, y))),\ (x, y) \in E_{ij}.$$

Therefore, for all $(x, y) \in E$, HVBRFIF $f_1$ satisfies

$$f_1(x, y) = s_{ij}(x, y) f_1(\vec{L}_{ij}^{-1}(x, y)) + s'_{ij}(x, y) f_2(\vec{L}_{ij}^{-1}(x, y)) + q_{ij}(\vec{L}_{ij}^{-1}(x, y)).$$

Moreover, $f_2$ satisfies $f_2(x, y) = \tilde{s}_{ij}(x, y) f_1(\vec{L}_{ij}^{-1}(x, y)) + \tilde{s}'_{ij}(x, y) f_2(\vec{L}_{ij}^{-1}(x, y)) + \tilde{q}_{ij}(\vec{L}_{ij}^{-1}(x, y))$.

### 5.2. Box-counting dimension of HVBRFIFs

In this section, we get a lower and upper bounds for the box-counting dimension of the graph of the HVBRFIF in the case where the data set be $P_0 = \left\{\left(x_0 + \frac{x_n - x_0}{n} i,\ y_0 + \frac{y_m - y_0}{m} j,\ z_{ij}\right) \in \mathrm{R}^3; i = 0, 1, \cdots, n,\ j = 0, 1, \cdots, m\right\}$, the extended data set $P = \left\{\left(x_0 + \frac{x_n - x_0}{n} i,\ y_0 + \frac{y_m - y_0}{m} j,\ z_{ij}, t_{ij}\right) \in \mathrm{R}^4;\ i = 0, 1, \cdots, n,\ j = 0, 1, \cdots, m\right\}$ with $\frac{|I_x|}{n} = \frac{|I_y|}{m}$, $x_{e(k)} - x_{s(k)} = y_{e(k)} - y_{s(k)} = \frac{\mu |I_x|}{n}\left(= \frac{\mu |I_y|}{m}\right)$ and $s_{ij}(x, y) s'_{ij}(x, y) \ge 0$, $\tilde{s}_{ij}(x, y) \tilde{s}'_{ij}(x, y) \ge 0$, $(x, y) \in E,\ i = 0, 1, \cdots, n,\ j = 0, 1, \cdots, m$. The row-stochastic matrix $M$ in (1) is irreducible.

Let us denote as follows:

$$P_{0x_\alpha} = \left\{\left(x_0 + \frac{x_n - x_0}{n}\alpha,\ y_0 + \frac{y_m - y_0}{m} j,\ z_{\alpha j}\right) \in \mathrm{R}^3;\ j = 0, 1, \cdots, m\right\},\ \alpha = 0, 1, \ldots, n,$$
$$P_{0y_\beta} = \left\{\left(x_0 + \frac{x_n - x_0}{n} i,\ y_0 + \frac{y_m - y_0}{m}\beta,\ z_{i\beta}\right) \in \mathrm{R}^3; i = 0, 1, \cdots, n\right\},\ \beta = 0, 1, \ldots, m,$$

Let us denote as follows:

$$\overline{\Omega}_\kappa = \max_{(x, y) \in E_\kappa} \{|s_\kappa(x, y)|, |s'_\kappa(x, y)|\},\ \overline{\tilde{\Omega}}_\kappa = \max_{(x, y) \in E_\kappa} \{|\tilde{s}_\kappa(x, y)|, |\tilde{s}'_\kappa(x, y)|\},$$
$$\underline{\Omega}_\kappa = \min_{(x, y) \in E_\kappa} \{|s_\kappa(x, y)|, |s'_\kappa(x, y)|\},\ \underline{\tilde{\Omega}}_\kappa = \min_{(x, y) \in E_\kappa} \{|\tilde{s}_\kappa(x, y)|, |\tilde{s}'_\kappa(x, y)|\},$$
$$\kappa = \tau(i, j) = i + (j - 1)n,\quad i = 1, \cdots, n,\ j = 1, \cdots, m.$$

As the matrixes $(\underline{S} + \underline{\tilde{S}})C$ and $(\overline{S} + \overline{\tilde{S}})C$ are positive irreducible, by Perron-Frobenius Theorem, they have spectral radius $\underline{\lambda}$ and $\overline{\lambda}$ with $\underline{\lambda} \le \overline{\lambda}$, respectively, where $C$ is the connection matrix of (22) and

$$\underline{S} = \mathrm{diag}(\underline{\Omega}_1,\ \underline{\Omega}_2,\ \cdots,\ \underline{\Omega}_N),\ \underline{\tilde{S}} = \mathrm{diag}(\underline{\tilde{\Omega}}_1,\ \underline{\tilde{\Omega}}_2,\ \cdots,\ \underline{\tilde{\Omega}}_N),$$

$$\overline{S} = \mathrm{diag}(\overline{\Omega}_{11},\ \overline{\Omega}_2,\ \cdots,\ \overline{\Omega}_N),\ \overline{\overline{S}} = \mathrm{diag}(\overline{\overline{\Omega}}_1,\ \overline{\overline{\Omega}}_2,\ \cdots,\ \overline{\overline{\Omega}}_N).$$

We have the following theorem.

**Theorem 6.** Let $f_1(x, y)$ be the HVBRFIFs constructed above. Suppose that there exist three interpolation points $(x_\alpha, y_{j_1}, z_{\alpha j_1})$, $(x_\alpha, y_{j_2}, z_{\alpha j_2})$, $(x_\alpha, y_{j_3}, z_{\alpha j_3}) \in P_{0x_\alpha}$ $(y_{j_1} < y_{j_2} < y_{j_3})$ (or $(x_{i_1}, y_\beta, z_{i_1\beta})$, $(x_{i_2}, y_\beta, z_{i_2\beta})$, $(x_{i_3}, y_\beta, z_{i_3\beta}) \in P_{0y_\beta}$ $(x_{i_1} < x_{i_2} < x_{i_3})$) which are not collinear and that take $t_{\alpha j_1}$, $t_{\alpha j_2}$ and $t_{\alpha j_3}$ (or $t_{i_1\beta}$, $t_{i_2\beta}$ and $t_{i_3\beta}$) such that $(z_{\alpha j_k} - z_{\alpha j_l})$ $(t_{\alpha j_k} - t_{\alpha j_l}) > 0$ (or $(z_{i_k\beta} - z_{i_l\beta})$ $(t_{i_k\beta} - t_{i_l\beta}) > 0$ ), $k, l = 1, 2, 3, k \neq l$ and three points $(x_\alpha, y_{j_1}, t_{\alpha j_1})$, $(x_\alpha, y_{j_2}, t_{\alpha j_2})$, $(x_\alpha, y_{j_3}, t_{\alpha j_3})$ (or $(x_{i_1}, y_\beta, t_{i_1\beta})$, $(x_{i_2}, y_\beta, t_{i_2\beta})$, $(x_{i_3}, y_\beta, t_{i_3\beta})$) are not collinear. Then the box-counting dimension of the graph of $f_1(x, y)$ is as follows:

(*a*) *If* $\underline{\lambda} > \mu$, *then* $1 + \log_\mu \underline{\lambda} \le \dim_B Gr(f_1) \le 1 + \log_\mu \overline{\lambda}$,

(*b*) *If* $\overline{\lambda} \le \mu$, *then* $\dim_B Gr(f_1) = 2$.

## 6. Conclusion

The HVRFIF is more general than the HVFIF. In particular, the HVRFIFs with function contractivity factors have more flexibility than ones with constant contractivity factors. First of all we estimate errors of HVRFIFs on perturbation of function contractivity factors, which ensures the stability of the constructed HVRFIFs according to small change of function contractivity factors. Next, we get the upper and lower bounds of box-counting dimension of the HVRFIF and HVBRFIF under some conditions, which gives fractal structure of the HVRFIFs and HVBRFIFs. These results become theoretical basis of practical application such as computer graphics, approximation theory, image process, data fitting and so on.